\documentstyle{article}
\begin{document}
\bigskip

\begin{center}
{ \Large \bf On Two Moduli Problems Concerning Number
 of Points and Equidistribution over Prime Finite Fields
    }
\end{center}
\smallskip
\begin{center}
{\ Nikolaj M. Glazunov}  \end{center}
\smallskip
\begin{center}
{\rm Glushkov Institute of Cybernetics NAS \\
   03680 Ukraine Kiev-680 Glushkov prospect 40 \\
   Email:} {\it glanm@d105.icyb.kiev.ua }
 \end{center} \smallskip

   \begin{center} {\bf Introduction} \end{center}
   The subject matter of this communication lies in the area between
 moduli theory~\cite{HM:MC}  and arithmetic geometry over finite fields.
Let ${\bf B}$ be a class of objects. In our case these are
classes of  hyperelliptic curves of genus ${\it g}$ over prime finite
field ${\bf F}_{p} $
 and Kloosterman sums $ T_{p}(c,d),  c,d \in {\bf F}_{p}^{\ast}  $. 
Let $ S $  be a scheme. A family of objects parametrized by the $ S $ 
is the set of objects

    $$ X_{s}: s \in S, X_{s} \in {\bf B} $$

  equipped with an additional structure compatible with the structure
  of the base $ S $.  \\
  We shall consider two problems: \\
   (i) existence of precise (exact) bound
  for families of hyperelliptic curves over $ {\bf F}_{p} $; and \\
  (ii) equidistribution of angles of Kloosterman sums. \\
  \begin{center} {\bf Moduli and estimates for hyperelliptic curves
  of genus ${\it g} \geq 2$ over ${\bf F}_{p}$ } \end{center}

  Let

  $$ C: y^{2} = f(x) $$

 be an algebraic curve and let $ Disk(C) $ be the discriminant of $ f(x) $.
 Let $p \geq 3$ be a prime.
 Consider  hyperelliptic curve of genus ${\it g} \geq 2 $ over prime finite
 field $ {\bf F}_{p}$

   $$ C_{g}: y^2 = f(x),\  D(f) \neq 0. $$

 For projective closure of $ C_{g} $ the quasiprojective variety

 $$ S_{g,p} = \{ {\bf P}^{2g+2}({\bf F}_{p}) \setminus (Disk(C_{g}) = 0) \} $$

parametrizes all hyperelliptic curves of genus ${\it g} $ over ${\bf F}_{p}.$
By well  known Weil bound (affine case)

 $$  |\#C_{g}({\bf F}_{p}) -  p| \leq 2{\it g}\sqrt{p}. $$

where $\#C $ is the number of points on the curve $C$ over ground field.
As we can see from Weil (and some more strong) bounds, for $p \geq 17$ any
hyperelliptic curve of genus  ${\it g} = 2 $ has points in ${\bf F}_{p} $ for
these prime $p. $ Also for ${\it g} = 3$ every hyperelliptic (h) curve of
genus 3 has points in ${\bf F}_{p}$ for $p \geq 37.$ For
$p =  3, 5, 7, 11 $ there are examples of h-curves of genus $2$ that
have not points in ${\bf F}_{p}.$ By author's computations~\cite{GP:NP} ,
any h-curve of genus $2$ over ${\bf F}_{13} $ has
points in the field. Similarly, for $p =  3, 5, 7, 11, 13, 17, 19, 23 $
there are examples of h-curves of genus $3 $ that have not points
in ${\bf F}_{p} $~\cite{GP:NP,Z01:L}.    \\

 \begin{center} {\bf Problem of precise bound} \end{center}

 Let $\{{\cal{F}}_{g,p}\}$ be a family of moduli spaces which is
parametrized by parameters ${\it g}$ and $p$. Let
$c \in {\cal{F}}_{g,p} $ be an element and
let $\#c$ be a numeric characteristic of $c$. Let ${\bf b}(g,p,\#c)$
be a bound that is satisfied for all $c \in {\cal{F}}_{g,p} $.
Let ${\cal{P}}(c,{\bf b}(g,p,\#c))$ be a predicate on elements of
${\cal{F}}_{g,p}.$ \\
    We shall say that the bound ${\bf b}(g,p,\#c) $
{\it precisely (exactly) divides} family ${\cal F}_{g,p}$
if for any given ${\it g}$ there exists $p = p_{0}(g) $
such that for any  $p \geq p_{0}(g)$ and for all
$c \in {\cal F}_{g,p} $ the predicate
${\cal{P}}(c,{\bf b}) = TRUE,$ and for $p \leq p_{0}(g) $ exists
$c \in {\cal{F}}(g,p) $ with ${\cal{P}}(c,{\bf b}) = FALSE.$ \\
Let ${\cal{F}} = S_{g,p} $ be the quasiprojective variety
of hyperelliptic curves of genus ${\it g},
\; C_{g} \in S_{g,p}, \; \#C_{g}$ be the number of points
of $C_{g}$ in ${\bf F}_{p}.$ \\
 \\
 {\bf Problem} {\it Does the precise bound for family} ${\cal{S}}_{g,p}$
 {\it exists?  If the precise bound exists what is its representation?.}
 \\
\\
More elaborately we have the following situation:
let $$f(x) = x^{2n+1} + a_1 x^{2n} + \cdots + a_{2n} x + a_{2n+1}.$$ \\
Let $deg f = 2n + 1 (n = 1,2,3,...)$ or genus(C) $= g (g = 1,2,3,...)$.
 Below in examples the Mit'kin (M) bound is used. \\
If $deg f = 3$ or $g = 1$ then $p_0 = 3$ (by M-bound every this curve 
has points in ${\bf F}_{p}$ for $p \geq 5$). \\
If $deg f = 5$ or $g = 2$ then $p_0 = 13$~\cite{Gl:M80}  (by M-bound 
every this curve has points in ${\bf F}_{p}$  for $p \geq 17$). \\
If $deg f = 7$ or $g = 3$ then $p_0 = 29$ (conjecture) (by M-bound every 
this curve has points in ${\bf F}_{p}$ for $p \geq 31$).   \\
If $deg f = 9$ or $g = 4$ then $p_0 = ?$ (by M-bound every this curve 
has points in ${\bf F}_{p}$ for $p \geq 53$).\\
Let $\#c$ be the number of points of an algebraic curve $c$ over prime 
finite field. Then the predicate is 
((For all $c \in S_{g,p})\& {\bf b}(g,p,\#c)  \Rightarrow (\#c > 0$)).

 \begin{center} {\bf Problem of distribution of Kloosterman sums} \end{center}

 Let

  $$ T_{p}(c,d) = \sum_{x=1}^{p-1} e^{2{\pi}i(\frac{cx + \frac{d}{x}}{p})} $$

 $$ 1 \leq c, d \leq p - 1; \; x, c, d \in {\bf F}_{p}^{\ast} $$
 be a Kloosterman sum. \\
 By A. Weil estimate

  $$ T_{p}(c,d) = 2 \, \sqrt{p} \cos \theta_{p}(c,d) $$

  There are possible two distributions  of angles $\theta_{p}(c,d)$
on semiinterval $ [0, \pi ):$ \\
\\
a) {\it $p$ is fixed and $c$ and $d $ varies over ${\bf F}_{p}^{\ast};$
what is the distribution of angles  $ \theta_{p}(c,d)$ as
$p \rightarrow \infty \; ;$ }  \\
 \\
 b){\it $ \; c$ and $d$ are fixed and $p$ varies over all primes not
dividing $c$ and} $d.$ \\
\\
For the case $a)$ N. Katz~\cite{K:GS},  A. Adolphson~\cite{A:E} and
 Chai \& Winnie Li~\cite{LC:C} proved that 
$\theta$ are distributed on $[0,\pi)$ with density
$\frac{2}{\pi} \sin^{2} t.$ \\
It is interesting to compare results of computer experiments in cases
$a)$ and $b).$ Such computations~\cite{G83:KS} and~\cite{G97:R}
demonstrated that though in case $b)$ equidistribution is possible
but results of computation shows not so good compatibility with
equidistribution as in (proved) case $a).$    \\

\begin{center} {\bf Acknowledgments} \end{center}

I would like to thank V.M. Sidelnikov, Winnie Li and Mike Zieve for  
questions, notes and information.

\end{document}